\documentclass[a4]{article}

\usepackage{natbib}
\usepackage{amsmath}
\usepackage{amssymb}
\usepackage{bibentry}
\usepackage{graphics}
\usepackage{authblk}
\usepackage{graphicx}

\newtheorem{theorem}{Theorem}[section]

\newtheorem{proposition}[theorem]{Proposition}
\newtheorem{defn}{Definition}[section]
\newtheorem{proof}{Proof}[section]

\newcommand{\V}{{\mathit{FP}}}
\newcommand{\T}{{\mathit{TP}}}
\newcommand{\E}{\mathbb E}

\newcommand{\Rj}{{\mathcal{R}}}

\newcommand{\fdr}{{\mbox{FDR}}}

\newcommand{\fdp}{{\mathit{FDP}}}
\newcommand{\paren}[1]{\left(#1\right)}
\newcommand{\ind}[1]{{\mathbf{1}\set{#1}}}
\newcommand{\set}[1]{\left\{#1\right\}}
\newcommand{\st}{\;|\;}
\newcommand{\pfdp}{\mathit{pFDP}}
\newcommand{\pfdr}{\mathrm{pFDR}}
\newcommand{\pvalues}{p_1,\ldots,p_m}

\newcommand{\pvec}{\mathbf{p}}
\newcommand{\sumtau}{\sum_{i=1}^{\tau}\frac{1}{i}}
\newcommand{\summ}{\sum_{i=1}^{m}\frac{1}{i}}

\begin{document}
\title{A comprehensive error rate for multiple testing}

\author[1,2,3]{Djalel-Eddine Meskaldji
\thanks{djalel.meskaldji@epfl.ch}}
\author[1,2]{Dimitri Van De Ville}
\author[4]{Jean-Philippe Thiran}
\author[3]{Stephan Morgenthaler}
\affil[1]{FSB/MATHAA, \'Ecole Polytechnique F\'{e}d\'{e}rale de Lausanne (EPFL), Lausanne, Switzerland.}
\affil[2]{Institute of Bioengineering, \'Ecole Polytechnique F\'ed\'erale de Lausanne (EPFL), Lausanne, Switzerland}
\affil[3]{Department of Radiology and Medical Informatics,University of Geneva, Geneva, Switzerland}
\affil[4]{LTS5/STI, \'Ecole Polytechnique F\'{e}d\'{e}rale de Lausanne (EPFL), Lausanne, Switzerland.}

\renewcommand\Authands{ and }

\maketitle

\begin{abstract} {The higher criticism of a family of tests starts with the
    individual uncorrected p-values of each test. It then requires a procedure
    for deciding whether the collection of p-values indicates the presence
    of a real effect and if possible selects the ones that deserve closer
    scrutiny. This paper investigates procedures in which the ordered
    p-values are compared to an arbitrary positive and non-decreasing threshold
    sequence.}
\end{abstract}

\section{Introduction}
\label{Motivation}

\begin{figure}

\centerline{\includegraphics[width=15cm]{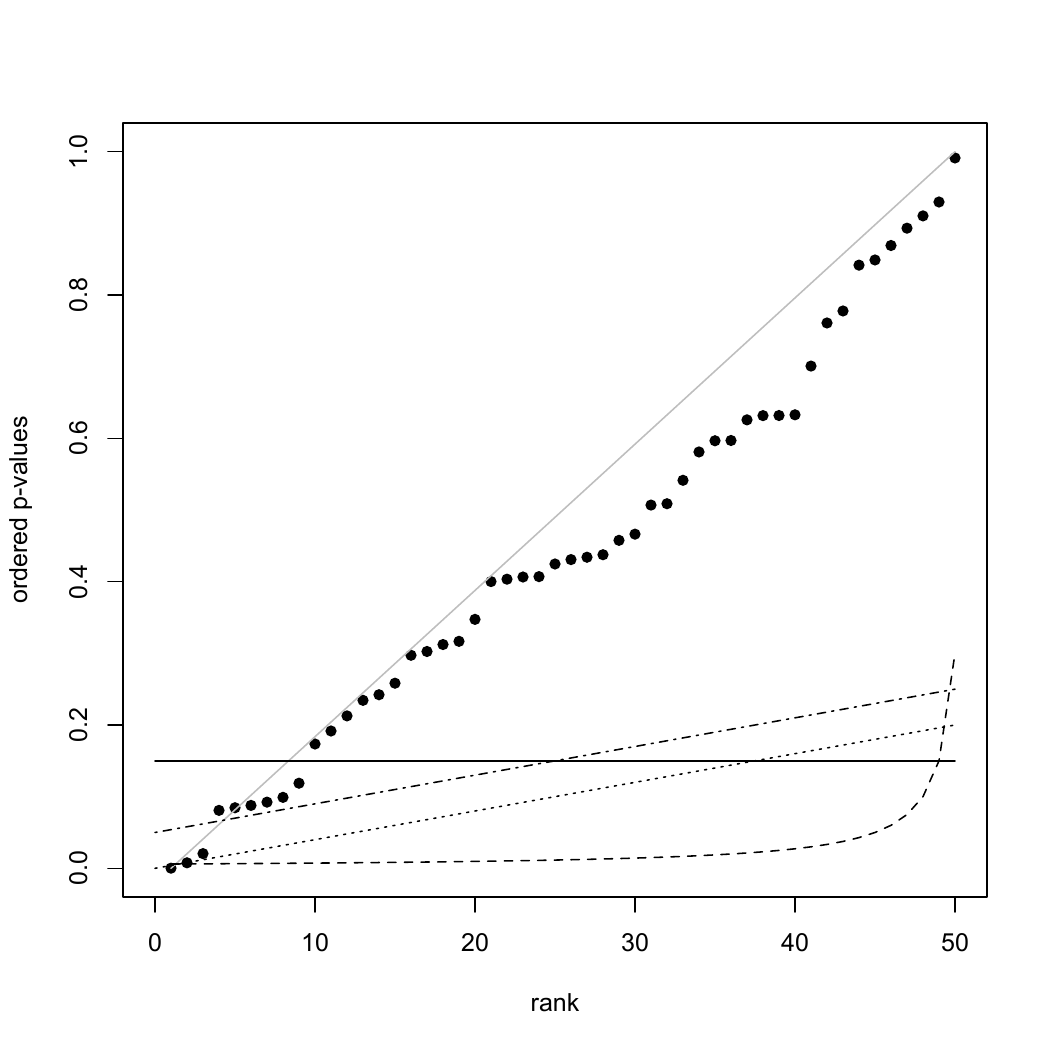}}

\caption{\label{thresholds} 
\small \em The plot shows for $m=50$ tests the ordered p-values
plotted against the rank. Besides the diagonal, three other curves
are shown, which have appeared in the literature as thresholds. They
are the constant threshold $t_r=0.15$ (solid line), the proportional
threshold $t_r=0.004 r$ (dotted line) and the curved threshold 
$t_r =0.3/(m+1-r)$ (dashed curve). With the
constant, 9 tests reject, with the proportional, 1 test rejects and with the
curved, 2 tests reject. In addition a linear threshold $t_r =0.05 +
0.004 r$ is shown.}

\end{figure}

When testing the existence of $m$ potential effects, it is customary to
compute an individual p-value $p_i$ for each null hypothesis $H_{0i},
\,i=1,\ldots,m$. A true effect tends to have a p-value close to zero, while
a null effect has a randomly assigned p-value between $0$ and $1$. Let
$m=m_0+m_1$ be the decomposition into the number of null effects and true
effects, respectively. The plot of the ordered p-values against the rank
provides a simple visual check to see whether any effect is present at
all. The resulting graph of $p_{(r)}$ vs $r$ for $r=1,\ldots,m$ carves out
an non-decreasing sequence of points in the rectangle $[1;m]\times[0;1]$ and
should, if no real effect exists, stay close to the diagonal from the lower
left corner to the upper right corner of the square. Procedures of higher
criticism can be based on a second curve or threshold sequence $t_r$, which stays
close to the lower edge of the square. If all p-values are above this
curve, none of the tests rejects. If any of the p-values are below the
threshold curve, there are several possible procedures. We will concentrate
on the maximum p-value below the curve. All the tests with p-value smaller
or equal to this maximum then reject. Figure \ref{thresholds} illustrates
this general idea. If such a procedure is applied, each effect we test is either
rejected or not rejected. Table \ref{general outcome} introduces the
contingencies that occur.

\begin{table}
\begin{center}
\begin{tabular}{c||ccc}
Hypotheses&not rejected& rejected&Total\\ 
\hline\hline
True& \textit{TN} & \textit{FP} & $m_{0}$ \\
False &\textit {FN} & \textit{TP} & $m_{1}$ \\
\hline
Total& $m-R$ & $R$ & $m$ \\
\end{tabular}
\end{center}
  \caption{\em General outcome of separating $m$ tests into non
    significant and significant ones. TP, FP denote the
    number of true, false rejections (true, false positives), respectively.  
    $R$ is the total number of rejections. TN, FN denote the number of
    true, false non-rejections (true, false negatives), respectively.}
\label{general outcome}
\end{table}

If we take the classical meaning of ``statistical significance at level
$\alpha$'' and reject all tests with $p_{(r)}\leq t_r=\alpha$, or
equivalently $p_i\leq \alpha$, the expected number of false rejections is
equal to $m_{0}\alpha$. If $m=1000$ and $\alpha=0.05$, we expect around 50
false rejections under the assumption that $m_1$ is moderate.
In Figure 1, the horizontal threshold has $\alpha=0.15$ and we expect
around 7.5 false rejections. The 9 rejections we actually found is close
enough to 7.5 to suspect that most or even all the rejections are false
positives. 

Even though this calculation makes it clear that performing many 
tests in parallel challenges the classical sense of ``statistically
significant'', a large number of claims are still published in the
scientific literature without a proper control \citep{BenjaminiSI2010}. To
make things worse, the above calculation neglects previously published or
unpublished studies that have tested the same or related hypotheses, thus
even further increasing the risk of spurious findings.

There is an obvious solution to the dilemma; we can simply replace
$t_r=\alpha$ by $t_r^{\text{Bonf}}=\alpha/m$. If we do this, the
expected number of false rejections is less than or equal to $m_{0}\alpha/m\leq
\alpha$. \cite{bonferroni1936} shows that in many circumstances this test
also controls another aspect, namely the probability of making one or more
false rejections. This stringent error control is one of the recommended
solutions to the dilemma of multiplicity, but can only be applied for
modest values of $m$, because otherwise the ability to find the real
effects of interest vanishes. 

More accommodating procedures were investigated. \cite{Holm1979}, for
example, showed that the curved threshold $t_r^{\text{Bonf}}\leq
t_r^{\text{Holm}}=\alpha/(m+1-r)$, lying along a convex function, provides
the same protection.\footnote{To be precise, Holm's procedure is a 
step-down procedure, that is , it rejects
$H_{(r)}$ as long as $p_{(r)}\leq t_r$ for $r=1,2,3,\ldots$ and stops
rejecting when for the first time the inequality fails.}
\cite{SIMES1986} went even further by considering the proportional choice
$t_r^{\text{Holm}}\leq t_r^{\text{Simes}}=r\alpha/m$. He showed that for
independent uniformly distributed p-values this still controls the
probability of one or more erroneous rejections. A breakthrough occurred
with \cite{BenjaminiFDR1995}. These authors recognized that it is useful to
study the proportional threshold in the context of another error metric,
namely the false discovery proportion $\fdp=\V/\max(R,1)$. They showed that
using the proportional threshold $t_r^{\text{Simes}}$ ensures that
$\mathbb{E}(\fdp)\leq \alpha$ (BH procedure).

In all these examples, the values of the thresholds lie between $\alpha/m$
and $\alpha$. This is reasonable since these two bounds represent the
cases of ``severe correction'' and ``no correction''. In this paper, we
explore what happens, if the threshold $t_r$ satisfies $\alpha/m \leq
t_r\leq \alpha$ and is non-decreasing,
for example, the linear threshold $t_r^{\text{linear}}=A + Br$, where $A,B
\geq 0$ and $A+B > 0$?

\section{Main results}
\label{Scaling functions section}

If we are willing to contemplate new error metrics, the answer to ``what
happens with arbitrary thresholds?'' has an easy answer.  Let the number of
rejections be $R$, the largest rank such that $p_{(R)} \leq t_R$ (see
Fig. \ref{thresholds}) and zero
if all the p-values are larger than the threshold. We will show that with
this threshold,
\begin{equation}
\frac{1}{m}\E\left(1_{\{R>0\}}\frac{\V}{t_R}\right)\leq 1 \,,\label{eq:error_controlA}
\end{equation}
which means that a ratio of the number of false positives divided by a
positive and non-decreasing function of the total number of rejections is
controlled. In order to relate this inequality to familiar results, we
consider thresholds of the form 
\begin{equation}
t_r=s_r\alpha/m \,,\label{eq:shapeandthreshold}
\end{equation} 
where the $0 \leq s_r \leq m$ is a non-decreasing sequence, which we will
call the shape of the threshold, and $\alpha>0$ is an arbitrary tuning
constant that controls the severity of the error control. The bigger
$\alpha$, the larger the threshold and the larger the number of
rejections. With this notation, inequality (\ref{eq:error_controlA}) becomes
\begin{equation}
\E\left(1_{\{R>0\}}\frac{\V}{s_R}\right)\leq \alpha \,.\label{eq:error_control}
\end{equation}
For the thresholds we have discussed in Section 1, the inequality
(\ref{eq:error_control}) yields the results in Table \ref{control}. 
\begin{center}
\begin{table}
\begin{center}
\begin{tabular}{c||c|c|c}
Name&threshold $t_{r}$&shape $s_r$ &control\\ 
\hline\hline
Bonferroni&$\alpha/m$ &$1$&
$\E(1_{\{R>0\}}\V)=\E(\V)\leq \alpha$ \\ \hline
Holm&$\alpha/(m+1-r)$&$\frac{m}{m+1-R}$&$\E(1_{\{R>0\}}\frac{m+1-R}{m}\V)\leq\alpha$\\ \hline
BH&$r\alpha/m$&$r$&$\E(1_{\{R>0\}}\frac{\V}{R})\leq\alpha$\\ \hline
Linear&$A + Br$&$\frac{m(A+Br)}{\alpha}$&$\E(1_{\{R>0\}}\frac{\V}{m(A+BR)})\leq \alpha$\\ \hline
\end{tabular}
\end{center}
\caption{\em The control exerted by using the thresholds discussed in 
  Section 1. The Bonferroni and BH case are well-known, the other two are
  new. The Holm threshold controls the expected number of erroneous 
  rejections, but multiplied by a weight smaller or equal to 1 that
  decreases linearly as $R$ increases. The linear and proportional
  thresholds also control such a weighted expectation, with decreasing
  weights depending on $R$. }
\label{control}
\end{table}
\end{center}
Inequality (\ref{eq:error_control}) gives us a clearer understanding of
procedures of higher criticism based on thresholding p-values of families
of tests. The threshold is meant to capture the p-values belonging to true
effects from above. In Fig. 1, such p-values cluster near the bottom left
corner, that is, they are small and have low rank. At least this is the
situation we have in mind, even though this will not happen if the true
effects are not large enough. If we move to the right along the curve of
the ordered p-values, the ones belonging to true effects will become rare
and the p-values sequence cuts across the threshold.

The proportional thresholds $t_r=r\alpha/m$ are ideal in situations where
the number of true effects is small and the effects are just barely
detectable. In other situations, the number $R$ of detected effects can be
large or moderately large and the proportional thresholds will then find a
sizable number of false positives and loose their effectiveness. This was
our motivation for studying other thresholds. We are, however, interested
in threshold shapes $s_r$ that grow linearly in $r$ for small values, but
sub-linearly for larger values $r$. The prototype is the trimmed threshold
$s_r^{trunc} = \min(\tau,r)$ for $1 \leq \tau \leq m$.

\begin{defn}
\label{def.thresh}
Consider an arbitrary non-decreasing threshold $ t_r = s_r\alpha/m$
$(r=1,\ldots,m)$ with non-decreasing shape function $1\leq s_r \leq
m$. The matched error metric is the penalized false discovery proportion
$$\mathit{pFDP} = \frac{\V}{s_R}\,, \text{ if } R\geq 1 \text{ and }
=0\,,\text{ otherwise}\,.$$
We call the expectation of the $\mathit{pFDP}$ the penalized false
discovery rate, $\pfdr=\E\paren{\pfdp}$. 
\end{defn}
In an older version of the manuscript \citep{MeskCER2011}, we show that with particular cases of the shape function, one can recover a large number of type I error rates that have been introduced in the literature. \\
Our results are valid for $m\geq 2$ and cover the $\pfdr$
control via step-up procedures\footnote{That is, procedures that reject all
  null hypothesis up to the last up-crossing of the threshold sequence.}
under independence, positive dependence and general dependence, as well as
the control via weighted procedures \citep{BH97weightedMTP,
  Genovese2006WFDR}. In a weighted procedure, an a priori non-negative
weight $w_{i}$ is assigned to each null hypothesis $H_{0i},$ and the raw
p-values $p_{i}$ are replaced by the weighted p-values $q_{i}=p_{i}/w_{i}$
if $w_{i} >0$ and $q_{i}=1$ if $w_{i} = 0$. A hypothesis of big weight is
more easily rejected than a hypothesis of small weight. The weights have to
be standardized so that $\sum_{i=1}^{m} w_{i}=m$, which includes $w_i\equiv
1$ as a possibility. The threshold can also be interpreted as a weight,
albeit one that applies to the rank of the hypothesis as determined by its
p-value. If we re-write $p_{(r)}\leq s_r \alpha/m$ as $p_{(r)}r/s_r\leq
r\alpha/m$, we see that the method acts as a weighted BH procedure with
weight proportional to $w_i\propto s_\text{rank of hypothesis i}/\text{rank
  of hypothesis i}$. Since these weights depend on the
p-values, the condition $\sum_{r}^{m}s_r/r=m$ does not guarantee the
control of the FDR.
\section{Discussion, Examples and Simulations}
\subsection{Relation to the literature}

Other studies have proposed threshold sequences different from the horizontal or
the proportional, among them \citet[p. 508 and
p. 513]{Genovese2002operating}, where it is shown that asymptotically the
proportional one  
is optimal in some sense, \citet{finner2009false}, who derive an
alternative threshold sequence based on asymptotic arguments, and
\citet{roquain2011exact}, who investigate the operating characteristics of the
FDR under an increasing threshold sequence. But none of these papers makes
the connection to a generalization of the error metric.  Metrics related
to ours were introduced in \cite{VDL2004} and described in \citet[p. 238
and ch. 6, 7]{DudoitLaan2007}. These references consider transformations
that involve both $\V$ and $R$, while we concentrate on transforming $R$
alone.

\subsection{The choice of the threshold sequence}

The freedom offered by using any threshold sequence is 
practical and useful. It gives users the
opportunity to choose on a very fine scale between tests that merely
provide a screening of a large numbers of potential 
effects and on the other hand tests that exert a strict control on the
number of false discoveries, but have difficulties in finding
moderately large true effects.

When applying a procedure that exerts FDR control, we can expect to have
$\V\approx \frac{m_{0}}{m}\alpha R$ for $m>>1$ \citep{roquain2011exact},
which means that the number of false effects ``discovered'' increases as a
fraction of the total number of effects found. The Bonferroni procedure
satisfies $\E(\V)\leq m \alpha/m = \alpha$, independently of $R$. Since
none of the error metrics is superior in all aspects, the users may wish to
have several distinct controls achieved by one procedure. On general
principles, an ideal procedure is as powerful as possible when it comes to
detecting true effects and, at the same time, controls the expected number
of false positives. The BH procedure does achieve this for
smallish $R$, but at the risk of large values of $\V$ when $R$ becomes
large. To avoid this, one can add strong control of $\E(\V)$ at a relaxed
level $\tau \alpha \; (\tau \geq 1)$ as \cite{HommelHoffmann88}
guarantees. In the context of penalized false discovery procedures, the
truncated linear threshold,
$t_r^{\text{trunc}}=s_r^{\text{trunc}}\alpha/m$, with shape function
$s_r^{\text{trunc}}=\min(\tau,r)$ provides a solution. The corresponding
step-up procedure with threshold sequence $ t_r^{\text{trunc}}$ controls
the FDR to be less than $\alpha$ and controls $\E(\V)$ to be less than
$\tau \alpha$. We might call this type of control the truncated FDR.  This
is an appropriate choice, if the number $m$ of hypotheses being tested is
very large. Applications include functional magnetic resonance images in
the neuroscience as well as many genomic or other ``omic'' studies
\citep{Meskaldji2015NI,Meskaldji2013NeuroImageReview}.

Note that the truncated FDR could be written as 
\begin{multline*}
\E \paren{\frac{\V}{ s_R^{\text{trunc}}}}=\E \paren{\frac{\V}{\min(\tau,R)}} =\max
\left[ \E \paren{\frac{\V}{R}}, \E \paren{\frac{\V}{\tau}}\right]\\
=\max\left[ \fdr, \frac{1}{{\tau}}\E \paren{\V}\right]\,.
\end{multline*}
A similar criterion (with smooth transition) is obtained by minimizing the
expected number of false positives while simultaneously controlling the false
discovery rate. This leads to a family of mixed error rates indexed by $0<\epsilon<1$ 
$$\text{MER}_{\epsilon, \tau}= \epsilon \E (\V /\tau) + (1-\epsilon)
\E(\V/R)\,,$$ 
where $\epsilon$ tunes the transition between the two controls. 
To derive a step-up procedure that controls the MER, we use Propositions \ref{PartI:
  SEV control indep} and \ref{PartI: SEV control pos. dep.}.  We
have 
$$\text{MER}_{\epsilon, \tau}= \E \left [\V \paren{\frac{R}{\epsilon R
      /\tau + (1-\epsilon) }}^{-1}\right]\,,$$ 
which we recognize as a pFDR with shape sequence
$$s_r^{\epsilon, \tau}=\frac{r}{\epsilon r /\tau + (1-\epsilon) }\,.$$
This sequence has a roughly linear rise for small $r$  and then flattens out to
reach the horizontal asymptote $r=\tau/\epsilon$.
Control at level $\alpha$, when $m$ hypotheses are tested, is obtained by the thresholds
$$t_{r}^{\epsilon, \tau}=s_r^{\epsilon, \tau}\alpha /m=\frac{r}{\epsilon r/\tau + 
(1-\epsilon) }\alpha/m\,.$$ 

\subsubsection{Simulations}
\begin{figure}[ht!] 
   \centering
   \includegraphics[width=\textwidth]{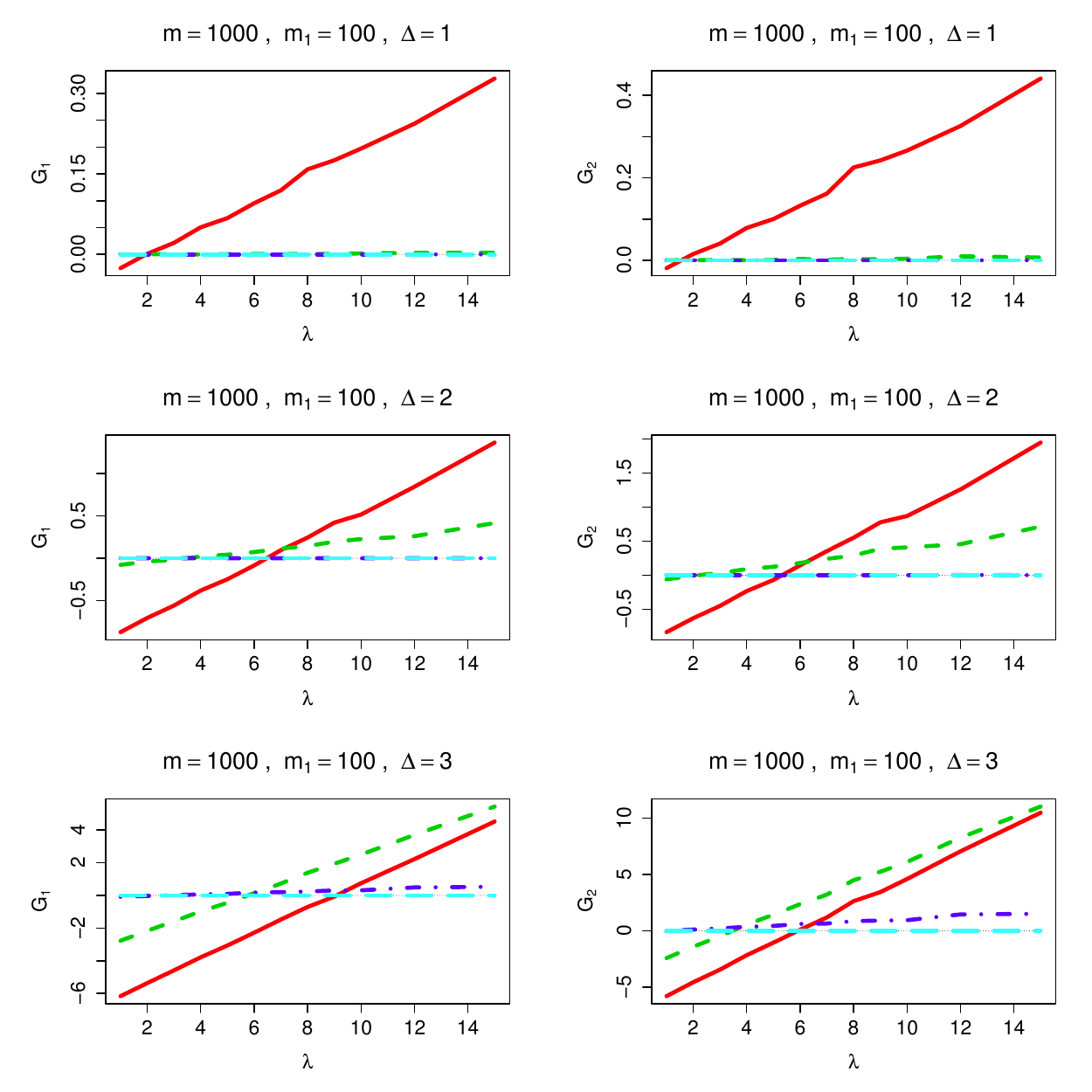}
   \caption{Behavior of the gain $\mathcal{G}_{1}$ and
   the gain $\mathcal{G}_{2}$ of different procedures compared to BH, as a function of the price $\lambda$. Here we show the difference between the gain function of each procedure and the one corresponding to the BH procedure. The number of tests is $m=1000$ and the number of false hypotheses is
     $m_{1}=100$. The effect size $\Delta$ is either 1,2, or 3. The curves
     show the simulated gains minus the gain of the BH procedure (zero line), Bonferroni (continuous), truncated linear with $\tau=5$
     (small dashes), $\tau=20$ (small-big dashes) and $\tau=100$ (large dashes).}
   \label{1000_100} 
\end{figure}
\begin{figure}[ht!] 
   \centering
   \includegraphics[width=\textwidth]{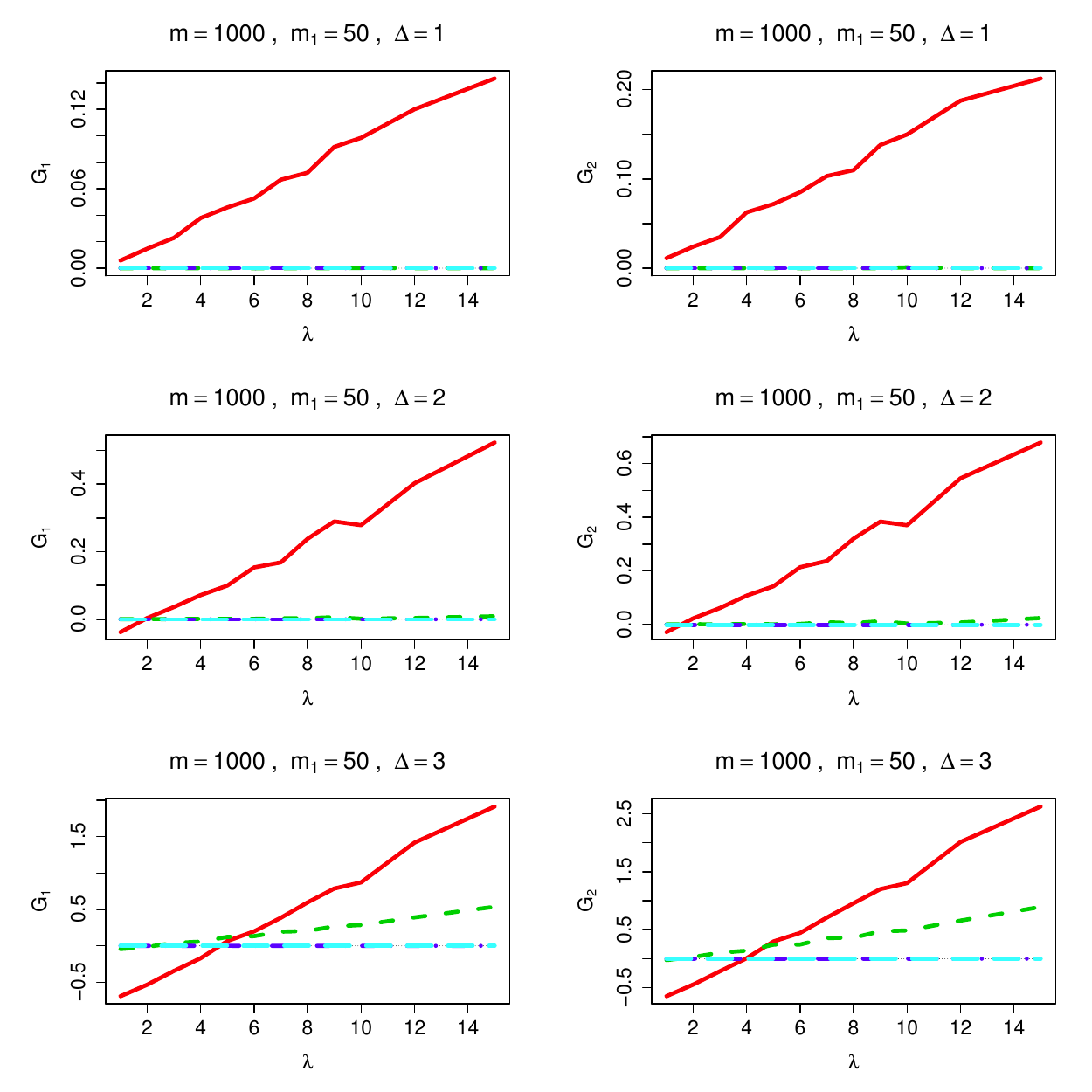}
   \caption{Behavior of the gain $\mathcal{G}_{1}$ and
   the gain $\mathcal{G}_{2}$ of different procedures compared to BH, as a function of the price $\lambda$. Here we show the difference between the gain function of each procedure and the one corresponding to the BH procedure. The number of tests is $m=1000$ and the number of false hypotheses is
     $m_{1}=50$. The effect size $\Delta$ is either 1,2, or 3. The curves
     show the simulated gains minus the gain of the BH procedure (zero line), Bonferroni (continuous), truncated linear with $\tau=5$
     (small dashes), $\tau=20$ (small-big dashes) and $\tau=100$ (large dashes).}
   \label{1000_50}
\end{figure}

\begin{figure}[ht!] 
   \centering
   \includegraphics[width=\textwidth]{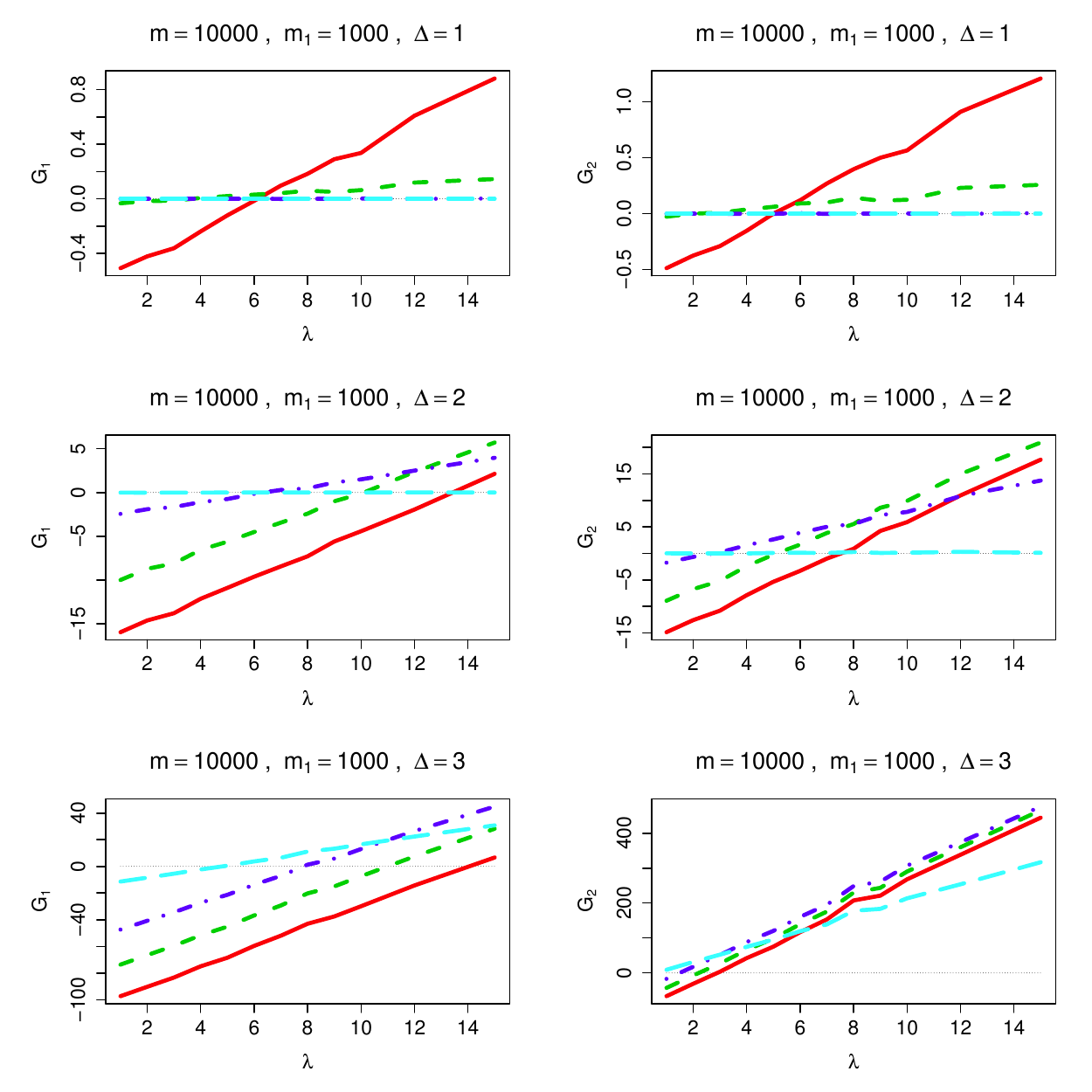}
   \caption{Behavior of the gain $\mathcal{G}_{1}$ and
   the gain $\mathcal{G}_{2}$ of different procedures compared to BH, as a function of the price $\lambda$. Here we show the difference between the gain function of each procedure and the one corresponding to the BH procedure. The number of tests is $m=10000$ and the number of false hypotheses is
     $m_{1}=1000$. The effect size $\Delta$ is either 1,2, or 3. The curves
     show the simulated gains minus the gain of the BH procedure (zero line), Bonferroni (continuous), truncated linear with $\tau=5$
     (small dashes), $\tau=20$ (small-big dashes) and $\tau=100$ (large dashes).}
   \label{10000_1000}
\end{figure}

\begin{figure}[ht!] 
   \centering
   \includegraphics[width=\textwidth]{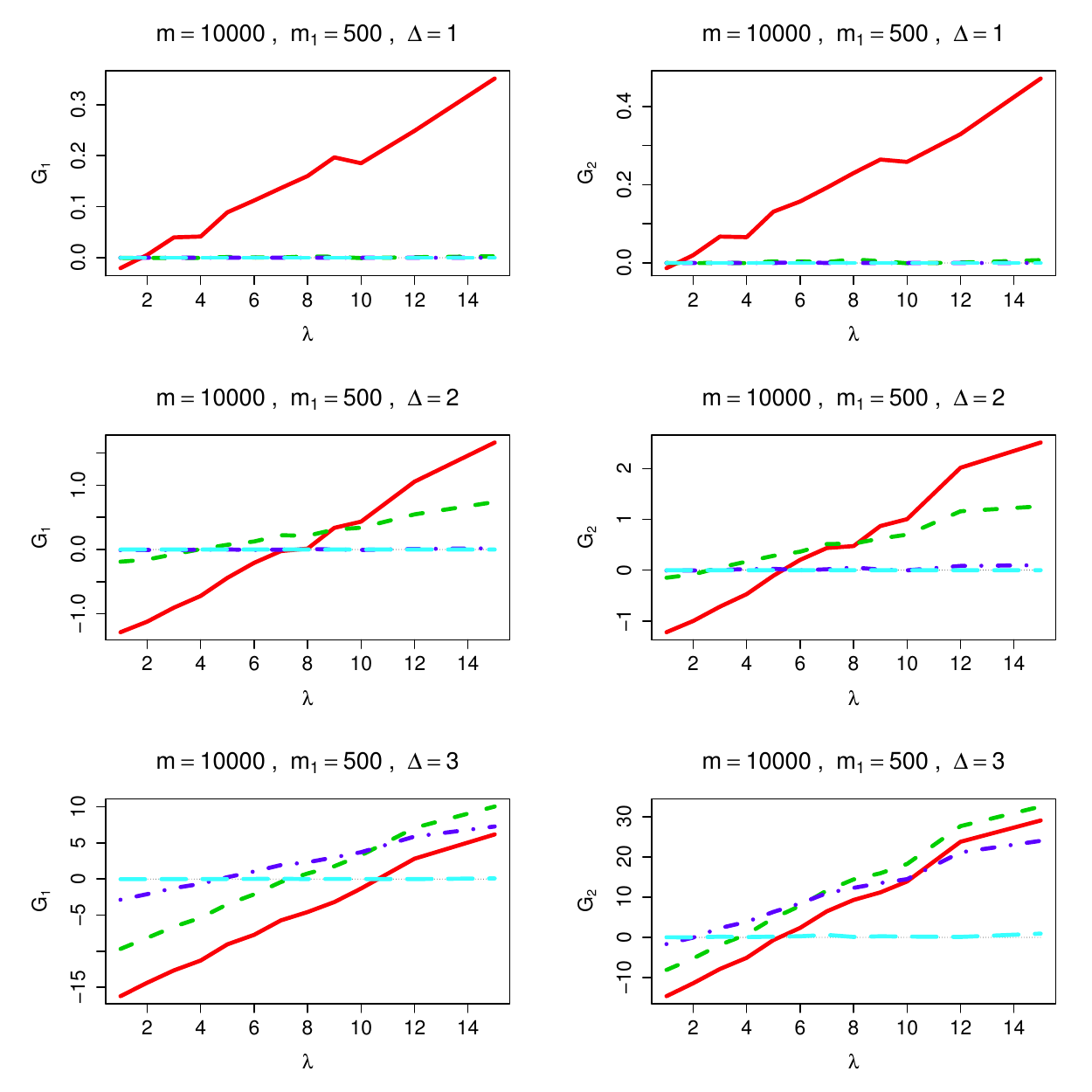}
   \caption{Behavior of the gain $\mathcal{G}_{1}$ and
   the gain $\mathcal{G}_{2}$ of different procedures compared to BH, as a function of the price $\lambda$. Here we show the difference between the gain function of each procedure and the one corresponding to the BH procedure. The number of tests is $m=10000$ and the number of false hypotheses is
     $m_{1}=500$. The effect size $\Delta$ is either 1,2, or 3. The curves
     show the simulated gains minus the gain of the BH procedure (zero line), Bonferroni (continuous), truncated linear with $\tau=5$
     (small dashes), $\tau=20$ (small-big dashes) and $\tau=100$ (large dashes).}
   \label{10000_500}
\end{figure}

To complete this section, we present some simulations that demonstrate the
advantage of using the truncated FDR. We compare three multiple testing
procedures: the Bonferroni procedure, the most stringent, the BH procedure,
the least stringent, and the truncated linear procedure which bridges the gap
between theses two extremes. We restrict our simulations to the case of
shifted Gaussians, that is, if the null hypothesis is true, the test
statistic follows a Gaussian law with expectation 0 and variance 1, while
false hypotheses are Gaussians with expectation $\Delta>0$ and variance
1. The values of the test statistics are independent. We will call $\Delta$
the effect. 

The general goal of any multiple testing procedure consists in making
$\T$ large while keeping $\V$ small. These two types of rejections are
opposites of each other, but asymmetrical opposites.  Among multiple
testing procedures that satisfy a given error control, the ones that
reject a maximal number of hypotheses are preferred. This is often the
only measure of optimality that is investigated. We consider instead a
simple monetary utility function with expected value
\begin{equation}
\mathcal{G}_{1}=\E(\T)-\lambda \E(\V)\,,
\end{equation}
in which each true discovery leads to a gain of 1 and each false discovery
to a loss of $\lambda \geq 1$. This approach might be unfamiliar to
statisticians, who are used to maximizing power under control of the false
rejections, but it can be shown that the two are closely related in the
Neyman-Pearson context. Our criterion allows a mixture of different error
rates and will pick the one best adapted to $\lambda$. In the philosophy of
multiple testing, $\lambda>1$, because the subsequent investigation of
any discovery is expensive and being on the wrong track may be a costly
mistake. 

It is stated in \cite{BenjaminiFDR1995} that the BH procedure is
optimal in the sense that under the FDR control, it maximizes the number of
rejections. In situations where $m_{1}$ and the effect are 
large, the number of true rejections is expected to be large. 
To maximize the number of rejections, counting both true and false,
the BH procedure can be expected to add rejections of null hypotheses,
even though the corresponding p-values are no longer small and we in
fact reject a sizable proportion of false ones. \\Another monetary
utility function to counteract such behavior has expected value
\begin{equation}
 \mathcal{G}_{2}=\E(\T)-\lambda \E(1+2+\ldots+\V)= 
\E(\T)-\frac{\lambda}{2} \paren{\E(\V)+\E(\V^{2})}\,,
 \end{equation}
where the loss due to a false discovery increases with each false
discovery. This gain function does not only consider the expected number of
false positive, but also its variance. 

We compare the gain functions of different procedures to the of the BH procedure, as a function of
the price $\lambda$. Figures 2, 3, 4 and 5 show the difference between the gain functions of the Bonferroni procedure, the truncated linear procedure with different values of $\tau$ (5, 20 and 100), and the BH gain functions, for various choices of the number of tests $m$, of the
number of false null hypotheses $m_1$ and of the effect $\Delta$. Any curve above zero indicates that the BH procedure is worse than the corresponding competing procedure, and vice versa for below zero. For both gains, the truncated linear procedures compensate a considerable part
of the power lost by the Bonferroni procedure without committing a large
number of false positives. Both the truncated linear and the Bonferroni
procedure are much more stable in terms of control of false
positives. The BH procedure seems to be appropriate only when the price
$\lambda$ of a false discovery is small. In some cases, the BH procedure
is optimal only when $\lambda =1$, which corresponds to the case where a
false rejection and a false acceptance of a null hypothesis are considered
to be of equal harm.  

As expected, the shortcomings of the BH are very pronounced for large
$m$. In the sparse case ($m_{1}=500$), the BH procedure is not much more
powerful than the other procedures when the raw effect is weak
($\Delta=1$). When the effect becomes moderate ($\Delta=2$) or strong
($\Delta=3$), the truncated procedures behave better than both the BH and the
Bonferroni procedure for reasonable values of the price $\lambda$. When the
number of alternative increases ($m_{1}=1000$ or $2000$), the performances
of the truncated methods are still good. However, the BH procedure shows a
marked performance instability, especially for $m_{1}=2000$ and
$\Delta=3$. \\
The simulations suggest the use of the
truncated linear with a value of $\tau$ that is situated between 20\% and
50\% of $\sqrt{m}$. 
\subsubsection{Examples}

\begin{figure}
\includegraphics[width=13cm]{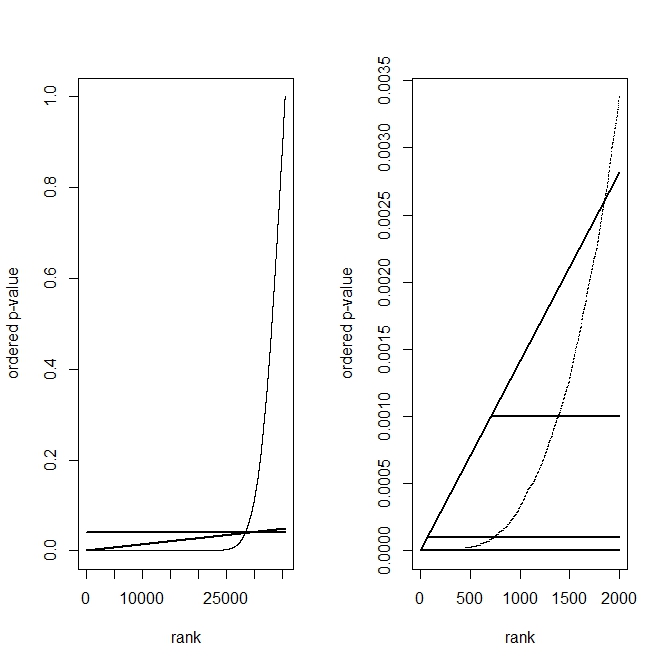}
\caption{Ordered p-values corresponding to gene expression level of $m=
  35512$ genes in two tissues of male mice, some fed a healthy diet, some a
  high-fat diet. In the left-hand panel the p-values are obtained from
  gene-by-gene comparison between muscle and liver, both diets
  confounded. The right-hand panel shows the p-values obtained when
 comparing the two diets in liver cells. }

\label{FigMice}
\end{figure}

Gene expression levels of $m=35512$ genes were measured in two tissues of
male mice who were fed either a normal diet or a high-fat diet
\citep{GeneDataCell2014}. We first consider an initial exploration of the
gene-by-gene comparison of the two tissues (strong difference). Two-sample
t-tests were performed to compare the expression levels in the muscle
versus the liver. The results are given as two-sided p-values, one for each
gene. The left-hand panel of Figure \ref{FigMice} shows that a large
majority of the genes exhibit a significant differential effect between the
two tissues. The Bonferroni procedure with $\alpha=0.05$ finds $21020$
differentially expressed genes, whereas the BH procedure with
$\alpha=0.05$ identifies $28436$ such genes. Under Bonferroni we stop
reporting positives early, in order to avoid the declaration of false
discoveries. Under BH, assuming that $m_1=29000$, one would expect
approximately $261$ falsely identified genes, namely $0.05\times28436\times
(6512/35512)$, whereas the Bonferroni procedures has an approximate
expectation of $0.05\times(6512/35512)=0.009$ false positives.

It is evident, that these two procedures cannot be compared
directly. If we were to match the approximate expected number of false
discoveries, the Bonferroni procedure would need $\alpha = 1422$, rather
than $\alpha=0.05$ and this procedure detects $28437$ significances, that
is, the two threshold sequences meet the ordered p-values almost at the
same point. The left-hand panel of Fig.~\ref{FigMice} shows these two threshold
sequences. The truncated procedure with $\tau=1422$, in which we decide to
apply FDR control up to a certain level and then switch to Bonferroni
control would, of course, give exactly the same result. This
illustrates a fact pointed out by  \cite{Genovese2002operating},
namely that the value where the
threshold sequence crosses the ordered p-values is the only aspect of
interest when investigating asymptotic error rates of the threshold
sequence.

Commenting on this first example, we can say that when many null hypotheses
are rejected, the BH procedure is a questionable choice for
performing multiple testing. By definition of the control metric, a
large number of rejections means that a large number of false rejections
are tolerated. The Bonferroni procedure with $\alpha=0.05$ tolerates none
or at most a few false discoveries, but can only achieve this feat by
limiting the reported number of significances. Among these first $21000$
positives the fraction of false ones is negligible. The BH
procedure with the same $\alpha=0.05$ has roughly an additional $7500$
rejections. Among these, we estimate that $261/7500=3.5\%$ are
false. The local rate of false discoveries increases when we only
consider the last 7500 rejections.

As a second example, the right-hand panel of Figure \ref{FigMice} shows the
(rounded) p-values one obtains when comparing the two diets in the
liver. In this case, we only take the data from one tissue and we thus have
only half the previous number of observations. Only a few of the $35512$
genes are differentially expressed, that is, the effect is weaker than in
the first example, and because of this, we only show the smallest $2000$
p-values. The superposed curves show the BH threshold sequence with
$\alpha=0.05$, the Bonferroni threshold with $\alpha=0.05$ and truncated
threshold sequences with $\tau=71$ ($\tau\approx
38\%$ of $\sqrt{m}$) and $\tau=710$ ($\tau\approx
380\%$ of $\sqrt{m}$).  The FDR control results in $1853$ discoveries,
whereas the Bonferroni procedure discovers $353$ genes. The truncated
thresholds result in $737$ and $1386$ rejections, respectively. This second
example shows, in a case where the effect is moderate, the advantage of
using the truncated threshold over the Bonferroni procedure. Up to a
certain number of discoveries, we accept the FDR control and then switch to a relaxed 
Bonferroni. If the effect becomes sparse and a truncated threshold is used,
the up-crossing of the ordered p-values may happen in the linearly
increasing part of the threshold sequence. In this case, the truncated
threshold behaves similarly to the FDR procedure and is able to detect
sparse signals, whereas the Bonferroni procedure will miss
them. At the same time, it offers better control. 

These examples illustrate how to use a truncated linear procedure in order to
increase the power considerably, but without incurring a punishingly high
expected number of false rejections. The optimal choice of the $\tau$
parameter depend on the price that the researcher is willing to pay when
committing a false discovery and on the (unknown) effect sizes.
\section{Results, proofs, and further discussions}

The theory developed in \cite{Blanchard2008TwoConditions}
can be adapted to our error rate. Three of their lemmas
form the basis, together with the concept of \emph{self consistency}. 

Any step-up procedure is such that the rejection set satisfies $\Rj = \{ i
\in I=\{1,\ldots,m\} : p_i \leq T \}$, where $T$ is an upper bound which
usually depends on the p-values $\pvalues$. To prove control, we need to bound
\begin{align}
  \nonumber \pfdr(T) =& \E\left[\frac{\V}{s_R} \ind{R>0}\right] \\=
  &\sum_{i \in I_{0}} \E\left[\frac{\ind{p_i \leq
        T}}{s_R}\ind{R>0}\right], \label{PartI: heuristic definition of
    the SEV}
\end{align}
where $I_0={\{i \in I \text{ such that } H_{i}\text{ is true}\}}$. The family of threshold
sequences $\Gamma: I \times \mathbb{R}^+ \mapsto \mathbb{R}^+$ of the form
$$\Gamma(i, r)=\frac{\alpha}{m} w_{i} s_r\,,$$
with $i$ denoting the null hypothesis and $r$ the rank of its p-value 
covers all cases of interest to us. It takes account of the weights
$\{w_i, i \in I\}$ and the shape $s_r$. 
Such a threshold sequence leads to a step-up procedure that satisfies 
$$ \Rj \subset \{i \in I \st p_i \leq \Gamma(i,R)\}\,,$$
where $\Rj$ is the set of rejections. This holds for all $R$ and is called the
self-consistency condition. 

\subsection{Control of the $\pfdr$ under independence}

\begin{proposition}\label{PartI: SEV control indep}
Assume that the $p$-values $\pvalues$ are independent random
variables and that non-negative weights $w_{1},\ldots,w_{m}$ are
given. It follows that the step-up procedure with threshold sequence
$t_{r}=s_r\alpha/m$ (see Def. \ref{def.thresh}) applied to
$q_{1}=p_{1}/w_{1},\ldots,q_{m}=p_{m}/w_{m}$ ($q_{i}=1$ if
$w_{i} = 0$), satisfies
$$\pfdr=\E(\V/s_R)\leq \alpha/m\sum_{i=1}^{m} w_{i}\,.$$
\end{proposition}

\begin{proof}
Denote by $\pvec_{-i}$ the set of $p$-values $\{p_j\st j\neq i\}$.
From (\ref{PartI: heuristic definition of the SEV}) and using the
self-consistency, we have
\begin{align*}
\pfdr (R) & = \sum_{i \in I_0}  \E\left[\frac{\ind{p_i/w_{i} \leq  s_R\alpha/m}}{s_R}\right] \\
& \leq\sum_{i\in I_0}  \E \left[ \E\left[ 
\frac{\ind{p_i\leq  \frac{\alpha}{m}
  w_{i}s_{R(\pvec)}}}{s_{R(\pvec)}} \middle\vert \pvec_{-i}
\right]\right]\\
& \leq  \frac{\alpha}{m} \sum_{i \in I_0} w_i\,,
\end{align*}
where we used the fact that given $\pvec_{-i}$, $p_i$ is
stochastically lower bounded by a uniform distribution in the
independent case. By case (i) of the Lemma 3.2 of
\cite{Blanchard2008TwoConditions} we obtain the result if we set
$U=p_i, \; g(U)=s(R(\pvec_{-i},U))$ and $c=\frac{\alpha}{m}
w_{i}$.
Note that $g(s_{R(\pvec_{-i},U)})$ is non increasing in $U$ since
$s_r$ is non-decreasing and $R$ is non-increasing in
$U$.
\end{proof}
\subsection{Control of the $\pfdr$ under positive dependence} 

\begin{proposition}\label{PartI: SEV control pos. dep.}
  Assume that the $p$-values $\pvalues$ are positively dependent
  (positively regression dependent on a subset as in 
  \cite{BY01}) random variables and that non negative weights
  $w_{1},\ldots,w_{m}$  
  are given. It follows that the step-up procedure with threshold sequence
  $ t_r=s_r\alpha/m$ (see Def. \ref{def.thresh}) applied to
  $q_{1}=p_{1}/w_{1},\ldots,q_{m}=p_{m}/w_{m}$ ($q_{i}=1$ if $w_{i} = 0$),
  satisfies
$$\pfdr=\E(\V/s_R)\leq \alpha/m \sum_{i=1}^{m} w_{i}\,, $$
where $R$ denotes the number of rejections.
\end{proposition}

\begin{proof}
By case (ii) of the Lemma 3.2 of \cite{Blanchard2008TwoConditions},  with $U=p_i$,
$V=s(R)$ and $c=\frac{\alpha}{m} w_{i} $, the result follows. 
\end{proof}

\subsection{Control of the $\pfdr$ under any dependence structure}

Under arbitrary dependence of the p-values, an additional correction to the
thresholds has to be performed. Let $\nu$ be a probability distribution on
$(0,\infty)$ and define for $r>0$
\begin{equation}
\xi(r) = \int_0 ^r u d\nu (u)\, . \label{equ_xi}
\end{equation}
We call the $\xi$ function \textit{the reshaping function}. It was
introduced in \cite{Blanchard2008TwoConditions}, who called $\xi(r)$
shape function, because they studied only the case of proportional
thresholds.

\begin{proposition}\label{PartI: SEV control general}
  Let $p_1,\ldots,p_m$ be the p-values computed for the $m$ null
  hypotheses. Non-negative weights $w_{1},\ldots,w_{m}$ are given. Let $\xi(r)$ be as described
  above. It follows that the step-up procedure with threshold sequence 
  $t_{r}=\xi(s_r) \alpha /m$ applied to
  $q_{1}=p_{1}/w_{1},\ldots,q_{m}=p_{m}/w_{m}$ ($q_{i}=1$ if $w_{i} = 0$),
  satisfies
  $$  \pfdr=\E(\V/s_R)\leq \sum_{i=1}^{m} w_{i}\alpha/m \,.$$
\end{proposition}

\begin{proof}
By case (iii) of the Lemma 3.2 of \cite{Blanchard2008TwoConditions} 
with $U=p_i$, $V=s(R)$ and $c= w_{i}\alpha/m$ the result follows.
\end{proof}

\subsection{Discussion on reshaping}

\begin{figure}
\centering\includegraphics[angle=270,width=\textwidth]{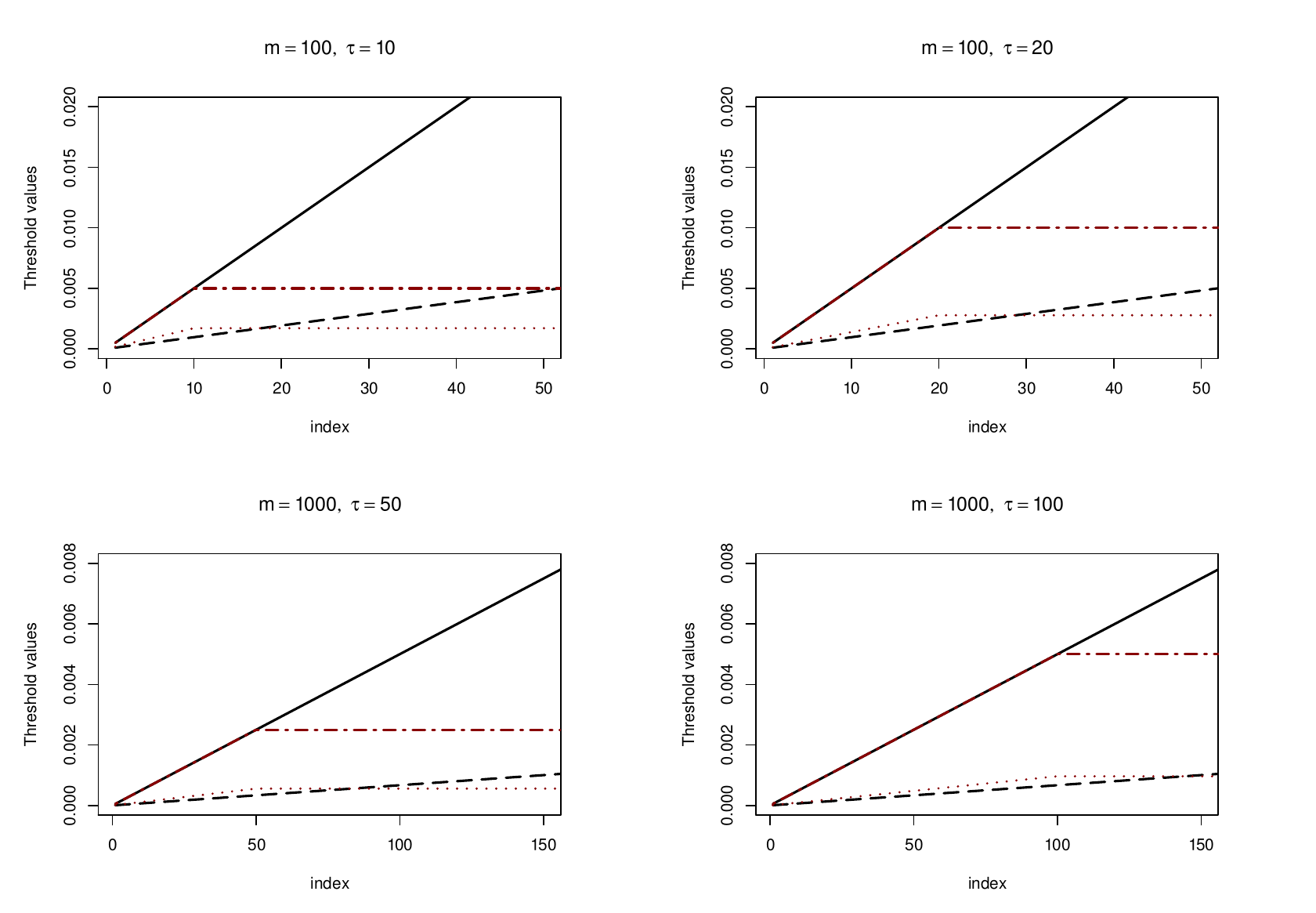}

\caption{\small \em Examples of threshold sequences for different values of
  $m$ and $\tau$. The BH procedure (continuous line) and the
  truncated linear procedure (big/small dashes), both for
  independent and positively dependent tests. The modified BH thresholds
  as described in \cite{BY01} (dashes) and the modified truncated linear
  thresholds as described in the Results section (points), both for the
  distribution free case.}
\label{thresholds2}

\end{figure}

The reshaping function $\xi(r)=r/(1+1/2+\ldots+1/m)$ is obtained for
the probability measure $\nu$ that puts point mass proportional to $1/r$ on
$r=1,\ldots,m$. The corresponding threshold sequence for the general
dependence case (Proposition \ref{PartI: SEV control
  general}) is
\begin{equation}
t_r=s_r\alpha/(m(1+1/2+\ldots+1/m))\,.\label{eq:BY}
\end{equation} 
The corresponding step-up testing procedure is a
penalized version of the procedure proposed in
\cite{BY01}. 
If we instead take $\nu$ as the probability measure with point mass
proportional to $1/1,\, 1/2,\ldots$ with support
$\{1,2,\ldots,\tau\}$, the reshaping function truncates, because
it is equal to 
$$\xi(r)=r\paren{\sumtau}^{-1}\text{ for }r\leq\tau\text{ and } 
\xi(r)=\xi(\tau)\text{ for }r>\tau\,.$$ 
The corresponding threshold sequence is 
$$s(r)^{\text{trunc}}\alpha\Big/\paren{m\sumtau}\,,$$ 
which exceeds (\ref{eq:BY}) for $r\leq{\lfloor}{\tau\summ/\sumtau
}{\rfloor}$.  This test thus increases the power in some situations even
though it exerts a provably stricter control (see Figure
\ref{thresholds2}). A particular case is $\tau=1$, where $\xi(r)\equiv
1$, and the threshold sequence is $\alpha/m$, that is, the Bonferroni
procedure.

\cite{Blanchard2008TwoConditions} propose other reshaping functions for the FDR case and comment on their use



\section{Further Results}
In addition to the freedom of choice of error rate, the duality that exists
between the error rate and the threshold
sequence $t_{r}=s_r\alpha /m$ together with the reshaping function $\xi(r)$,
brings new perspectives to the field of multiple testing procedures and
reveals new interesting results. In this section, we give some examples 
of corollaries and applications of our results. 

\subsection{Relation between error rates based on expectation \\ or on tail probabilities}

Define the penalized false exceedence rate ($p\text{FER}$) as 
$$p\text{FER}_{q}=\mathbb{P}\left(\frac{\text{FP}}{s_R}>q\right)\,.$$ 
By Markov's inequality
\begin{equation*}
p\text{FER}_{q}=\mathbb{P}\paren{\frac{\text{FP}}{s_R}>q}\leq\alpha 
\Rightarrow \mathbb{E}\paren{\frac{\text{FP}}{q \, s_R}}\leq \alpha\,,
\end{equation*}
from which it follows that the threshold sequence $qs_{r}\alpha /m$ provides
conservative control of
this error rate, under independence
(Theorem \ref{PartI: SEV control indep}) or positive dependence (Theorem
\ref{PartI: SEV control pos. dep.}). Under general dependency
structure, we replace $qs_R$ by $\xi(qs_R)$ (see 
Lemma \ref{PartI: SEV control general}). In an older version of the manuscript 
\citep{MeskCER2011}, we proposed step-down procedure that tightly controls the pFER, under different assumptions. 

An interesting choice is $q=0.5$, in which case the
$p\text{FER}_{q}$ becomes the median of $\mathit{FP}/s_R$ and we
note that in order to achieve control, we simply have to halve
$\alpha$. Consequently, if we wish to control the median of the FDP at
level $\alpha$, we can use the BH procedure with $\alpha/2$.

\subsection{Asymptotics}

Hochberg's procedure is a step-up procedure with Holm's threshold
sequence $t_r=\alpha/(m-r+1)$.
If we factor out an $\alpha/m$, the shape is seen to be  $s_r=m/(m-r+1)$. Our
results show that Hochberg's procedure controls 
$$\mathbb E \left( \frac{\mathit{FP}}{m} (m-R+1)\right)\approx
{\E(\V)}\text{ (for large }m)\,,$$
both under independence and positive dependence. 

Another example of this type is the asymptotic control of the optimal
rejection curve (ORC) introduced by \cite{finner2009false}, which controls
asymptotically the FDR. We give here a simple proof of that control. The
ORC threshold sequence is $t_r=r\alpha/(m-r(1-\alpha))$, which equals 1, when
$r=m$ and thus has $R=m$. With a slight modification it can be used in
practice. When
we factor out $\alpha/m$ we obtain $s_r=rm/(m-r(1-\alpha))$ and
we can conclude that the modified threshold controls the expectation  
$$\mathbb E \left(\frac{\mathit{FP}}{R} (1-R/m(1-{\alpha}))\right)\to
\text{FDR} \text{ as }m\to\infty\,. $$ 
This result is valid under independence and positive dependence.

\subsection{Assumption weakening}
For a step-up procedure associated with a specific threshold sequence, one
can find different error rates controlled by this procedure under different
assumptions on the dependence structure of the p-values. The
proportional threshold $t_r=r\alpha/m$, for example,  is obtained from $\alpha
\xi(s_r)/m$ when $\xi(r)$ is the
inverse of the shape sequence $s_r$, that is, $\xi(s_r)=r$. This shows
that the BH procedure controls the quantity
$$ \E\left[\frac{\V}{s_R}\right] =\E\left[\frac{\V}{\xi^{-1}(R)}\right]$$ 
at level $\alpha$ under any distribution of the p-values and for an
arbitrary choice of $\xi(r)$ or the shape sequence $s_r$. 
We are thus able to determine the error rate that is under control if
the independence or positive dependence assumption is
violated. This result is obvious for a linear $\xi(r)$
such as $\xi(r)=r/\summ$, in which case
$$\xi^{-1}(r) =r\summ=s_r\,.$$ 
The BH procedure, no matter the dependence structure of the p-values, controls
$$\E\paren{\V/\paren{R\summ}}\leq \alpha\,,$$ 
which immediately leads to 
$$\E\paren{\V/R}\leq \alpha \summ\approx\alpha\,\log (m)\,.$$ 
Our results thus inter-relate different error rates and conditions on the dependence structures of
p-values by the relation
$$ t_{r}=\xi(s_r)\alpha/m \equiv \tilde s_r\alpha/m, \text{for } r=1,\ldots,m.$$


\section*{Acknowledgments}
This work was supported by the Swiss National Science Foundation [144467, PP00P2-146318].\\
A part of this work was done while the first author was a Ph.D student
in the Signal Processing Laboratory (LTS5),  Ecole Polytechnique
F\'{e}d\'{e}rale de Lausanne (EPFL), Lausanne, Switzerland \citep{Meskaldji2013thesis}.\\ 
The data is from the Nestl\'e chair in Energy Metabolism (Prof. J. Auwerx)
and was collected by a graduate student (E. G. Williams).\\ 
The authors would like to thank Etienne Roquain for interesting comments
and suggestions.  
\bibliographystyle{chicago} 
\bibliography{Bib}

\end{document}